\documentclass[11pt]{article}
\usepackage{amsbsy,amsfonts,amsmath,amssymb,amscd,amsthm}
\usepackage{mathrsfs}
\usepackage[a4paper,text={6.1in,8in},centering,includefoot,foot=1in]{geometry}
\usepackage{amsmath,amsthm,amsfonts,amssymb,amscd}
\usepackage{graphics}
\usepackage{color}
\usepackage[colorlinks=true,citecolor=blue]{hyperref}
\usepackage{latexsym}
\usepackage{cite}
\usepackage{url}
\usepackage[dvips]{graphicx}

\small\normalsize
\theoremstyle{plain}
\newtheorem{theorem}{Theorem}

\newtheorem{corollary}{Corollary}

\newtheorem{definition}{Definition}

\newtheorem{lemma}{Lemma}

\newtheorem{proposition}{Proposition}

\numberwithin{equation}{section}

\begin{document}

\title{\textbf{Topological  structure of  the \\ sum  of  two  affine Cantor sets}}
\date{}
\author{ Mehdi Pourbarat\\
\\\small \emph{ Shahid Beheshti Uni., Tehran,  Iran.}} \maketitle
\maketitle
\centerline{\it Dedicated to my wife, \textbf{Mina}.}
\vspace{.5cm}

\hspace*{-.5cm}\textbf{Abstract}: 
We introduce a dense subset of affine Cantor sets, termed "generalized homogeneous Cantor sets." We show that for any two members of this class, if their sum set is not a Cantor set,  it has a dense interior.
If one of them is an affine Cantor set with two mappings, the same result holds.
In the contex of affine Cantor sets defined by increasing maps, we introduce a dense subset of their pairs, denoted by $\cal{D}$, such that for every $(K, K') \in \cal{D}$, there are five possible structures for their sum set:  a bilateral gap interval, an L, R, M-Cantorval, or a Cantor set.   
Finally,   we  present new pairs of  affine Cantor sets  which have stable intersection, while   do not satisfy the  Generalized Thickness Test. 
\\\\ 
\textbf{Keywords}:     Affine   Cantor sets,  Palis Conjecture, Newhouse's thickness.
\\AMS Classification: 28A78, 58F14.

\section{Introduction}\label{sec1}
Affine Cantor sets  have been regarded, on one hand as   a   natural generalization of the 
 homogeneous  Cantor sets   and   on the other hand as  a particular case of regular Cantor sets. An affine   Cantor set can be obtained via an iteration process on the finite union of many discrete   closed intervals  not necessarily same lengths  \cite{Fa}. In one dimensional  dynamical systems theory,  regular Cantor sets  may  appear   as repellers
of expanding  maps, for definitions and background, see \cite{Fa1}.  In the study of  two-dimensional diffeomorphisms, they also may  appear from intersections  of hyperbolic sets  with stable or unstable manifolds of its points   \cite{PT}. 
 The arithmetic   sum  of two regular Cantor sets  sometimes   play a key role  in various settings such as   dynamical systems, number theory,  real analysis,  spectral theory,  harmonic analysis (for each, see a number of references in \cite{P4}) and
the existence of patterns in fractals \cite{Ma,Y}.

In dynamical systems theory,  a  Cantor set $K$  is \textit{regular} or \textit{dynamically defined} if:
   \begin{itemize}
   	\item [(i)] there are disjoint compact intervals $I_{0}, I_{1},\ldots,I_{M}$ such that
   	$K\subset \bigcup_{i=0}^M I_{i}$ and the boundary of each $I_{i}$ is
   	contained in $K$,
   	\item [(ii)] there is a $C^{1+\epsilon}$ expanding map $\psi$
   	defined in a neighborhood of set $ \bigcup_{i=0}^M I_{i}$ such that
   	$\psi(I_{i})$ is the convex hull of a finite union of some intervals $I_{j}$
   	satisfying:
   	\begin{itemize}
   		\item [$\bullet$]  for each $0\leq i\leq M$ and $n$ sufficiently
   		large, $\psi^{n}(K\cap I_{i})=K,$
   		\item [$\bullet$]$K=\bigcap_{n=0}^{\infty}\psi
   		^{-n}( \bigcup_{i=0}^M I_{i})$.
   	\end{itemize}
   \end{itemize}
   The set $\{I_{0},I_{1},\cdots,I_{M}\}$ is, by definition, a Markov partition for $K$, 
   the set $D:=\bigcup_{i=0}^{M} I_{i}$ is  Markov domain   of $ K$ and $K_n=\psi
   		^{-n}( \bigcup_{i=0}^M I_{i})$ is called  the $n$-th step of the construction  of $K$.  A regular Cantor set $K$ is \textit{affine} if $D\psi$ is constant  on every interval $I_i$. It is called \textit{homogeneous}, if  $D\psi$ is constant  on the Markov domain $D$.   
   		
  Suppose that $K$ is an affine Cantor set  with  Markov domain  $\{I_{0},I_{1},\cdots,I_{M}\}$
 and  expanding map $\psi$ defined by $\psi_{|I_i}=:p_ix+e_i$ and also $K'$ is an affine Cantor set  with  Markov domain  $\{{I'}_{0},{I'}_{1},\cdots,{I'}_{N}\}$
\label{sec2}  and  expanding map $\psi'$ defined by $\psi'_{|I'_j}=:q_jx+f_j$.   Let $\Sigma:= \{0,1,2, \dots, M\}$ and $\Sigma':=\{0,1,2, \dots, N\}$. Let ${\cal{A}}$  be   the space of all  affine  Cantor sets $K_{p_0,\dots ,p_M}^a$as above. For given   affine  Cantor sets $K_1$ and $K_2$ with arbitrary convex hulls,  there is the pair  $(K,K')\in {\cal{A}}\times {\cal{A}}$  such that   $K_1-K_2$ is  an   affine translation of  $K-  K'$.  Thus   the space of all pairs affine  Cantor sets  is homeomorphic with $ {\cal{A}}\times {\cal{A}}\times \mathbb R\times \mathbb R$. We consider the Lebesgue measure induced on it by  identifying   $K \leftrightarrow(p_0,\dots ,p_M,a)\in\mathbb R^{M+2}$. 
   		
Suppose that $K$ and $K'$ are two regular   Cantor sets and $\lambda$ is a real number.  Define
 $$K -\lambda  K' :=\big\lbrace x-\lambda y\mid ~x \in K,~y \in K'\big\rbrace.$$
 Note that  $-K'$ is a regular Cantor set too, thus   the  sum set  $K +\lambda  K'$ can be considered as a  difference set \cite{P11}.  
The reason that we are interested  in considering  difference sets in the current form is that for each element $t\in K -\lambda  K'$ and $\lambda\neq0$,  the sliding image $ t +\lambda  K'$  of  $\lambda  K'$ intersects $K$, and  from this point of view investigating   the geometric structure of the set $K -\lambda K'$  is useful in studying the  bifurcation phenomena in  dynamical systems \cite {PT}. 

In  geometric measure theory, regular Cantor sets may be   defined as  the attractors of the  iterated function systems  (IFSs) which  satisfy the open set condition.  Assume that affine Cantor sets  $K$ and $K'$ are the attractors of the  IFSs $\{F_i(x)=r_ix+t_i\}_{i=1}^{M}$ and $\{F_j'(x)=r_j'x+{t'}_j\}_{j=1}^{N}$ with  convex hulls  $[0,~a]$  and  $[0,~b]$, respectively.  Throughout the paper,  abbreviation  $BD(A),~HD(A),~H^s(A)$,  $|A|$ and $A^\circ $ denote  the  Box dimension,  Hausdorff  dimension, $s$-dimensional  Hausdorff measure, Lebesgue  measure and interior  of $A$, respectively  \cite{Fa1}.
Cantorval sets were first introduced by  by Mendes and Oliveira \cite{MO}.
\begin{definition} \label{Z} Suppose that $A$ is a   perfect  subset of $\mathbb{R}$  such that any   gap has an interval adjacent to its right (resp. left) and is accumulated on the left (resp. right) by infinitely many intervals and gaps. Then $A$ is called
 L-Cantorval (resp. R-Cantorval). It  is called  M-Cantorval, if  any   gap is accumulated on each side by infinitely many intervals and gaps.
\end{definition}
Thus, if $A$ is   L-Cantorval,  then $-A$ is   R-Cantorval and vice versa. To state our result, we present the following two definitions.
\begin{definition} \label{Z} Suppose that $A$ is a   perfect  subset of $\mathbb{R}$  such that any  bounded  gap has two intervals adjacent to its two side. Then $A$ is called  bilateral gap interval. 
\end{definition}
Thus,  if the set $A$ is  a  finite  union of closed intervals, then  $A$ is   bilateral gap interval, Also  $A$ is   bilateral gap interval,  if and only if $-A$ does.

\begin{definition} \label{Z} An affine Cantor set $K$ is called a generalized homogeneous Cantor set (GHC) if, for all pairs $(i,j) \in \Sigma \times \Sigma$, we have $\frac{\log p_i}{\log p_j} \in \mathbb{Q}$. This set is termed "normal" if $p_0 = p_M$ and, for all $i \in \Sigma$, the condition $\frac{\log p_i}{\log p_0} \in \mathbb{N}$ holds. We denote the set of GHCs by $\mathcal{G}$ and their normal subset by $\mathcal{N}$.
\end{definition}
In Section~\ref{sec2},  we show that  for each  the   pair $(K,K')\in {\cal{G}}\times{\cal{G}}$, the set $K-K'$   has dense interior, if it is not Cantor set. Furthermore, this statement holds if $K'$ is an affine Cantor set generated by two maps (Theorem \ref{A}). We  show that for any pair $(K, K') \in \mathcal{N} \times \mathcal{N}$, the set $K-K'$ is  generically (i.e., from both topological and measure-theoretic) a bilateral gap interval  or  an  L, R, M-Cantorval or a Cantor set  (Theorem \ref{B}). Finally, we introduce a dense subset, such as $\cal D$, of the space $\{(K,K')\in {\cal{A}}\times{\cal{A}}|~ D\psi,D\psi'>0\}$ in which the set $K-K'$ lacks the five aforementioned properties   (Theorem \ref{BB}).

 We say that the Cantor set $K$ is close on the topology $C^{1+\epsilon}$  to a Cantor set $\widetilde{K}$ with the Markov partition $\{\widetilde{I_1},\widetilde{I_2},\cdots,\widetilde{I_M}\}$
defined by expanding map $\widetilde{\psi}$ if and only if, the extremes of $I_i$ are near the corresponding extremes
of $\widetilde{I_i},~ i=1,2,...,M$ and supposing $\psi\in C^{1+\epsilon}$ with Holder constant $C$, we must have $\widetilde{\psi}\in C^{1+\widetilde{\epsilon}}$ with Holder constant $\widetilde{C}$
such that $(\widetilde{C},\widetilde{\epsilon})$ is near $(C,\epsilon)$ and $\widetilde{\psi}$ is close to $\psi$
in the $C^1$ topology.

\begin{definition}  Regular Cantor sets $K$ and $K'$  (or pair $(K,K')$) have \emph{stable intersection} if for any pair of regular Cantor sets $(\widetilde{K} ,\widetilde{K'})$ near $(K ,K')$, we have
$\widetilde{K}\cap \widetilde{K'} \neq\emptyset$.
\end{definition} 
  
From \cite{HMP}, transferred renormalization operators or simply operators are  given by
\begin{equation}
(s,t)\overset{T_i}\longmapsto
(p_is,p_it-p_ia_i),\hspace{1cm}(s,t)\overset{T'_j}\longmapsto
(\frac{s}{q_j},t+a'_js).
\end{equation}  
The operators $T_i$ and $T'_j$ comutate together; $T_i\circ T'_j=T'_j\circ T_i$,  and also they have  inverses  
${T_i}^{-1}(s,t)=(\frac{s}{p_i}, \frac{t}{p_i}+a_i)$ and  ${T'_j}^{-1}(s,t)=(qs,t-qa'_js)$,  for  all $i\in\Sigma$  and $j\in\Sigma'$,

 In [MY], Moreira and  Yoccoz guarantee the existence of an open and dense subset  among  all   pairs of regular Cantor sets  in $C^{1+\epsilon}$ topology with the sum of their   Hausdorff dimension  bigger than one, such that their elements have  stable  intersection [MY]. For each pair of Cantor sets from this subset, the main challenge is to construct a recurrent compact set  of relative configurations of limited geometries corresponding to their renormalization operators. When Cantor sets are affine,  recurrent set  may state as below definition.
\begin{definition}  The compact set $\mathcal{R}\subset \mathbb{R^{*}}\times \mathbb { R} $ is a \emph{recurrent} set  for transferred renormalization operators if for every element of $\mathcal{R}$ there exist suitable composites of operators $(2.1)$, which transfers that element to $\mathcal{R}^{\circ}$.
\end{definition}
In fact, for every $(s, t)$ belonging to the interior of the  recurrent set $\mathcal R$, Cantor sets $K$ and $sK'+t$ have stable intersection. Moreover,   for each $(s,t)\in \mathcal R$,  $K-sK'$ contains an interval.
 We say that  the  pair $(s,t)$ is a \textit{difference pair }for  $K$ and $K'$,  if there exists a relatively compact  sequence $(s_{n}, t_{n})$  in  the space  $\mathbb{R}^* \times
\mathbb{R}$ with $(s_1,t_1)=(s,t)$  such that the point  $(s_{n+1},t_{n+1})$ is obtained from $(s_{n}, t_{n})$ by  applying one of  operators (1.1). Thus, if  $\mathcal {R}$ is a recurrent set, then each pair $(s,t)\in (s,t)\in \mathcal R$ is a  difference pair for  $K$ and $K'$. The following assertion  is proved in \cite{P1}:
\begin{itemize}
\item [$\bullet$] $t  \in K-sK'$   if and only if  $(s,t)$ is a difference pair for  $K$ and $K'$.
\end{itemize}

Suppose that $K$ is a regular Cantor set and   $U$ be a bounded gap  of  $K$.  Also, suppose that $L_U$ and  $R_U$ are the intervals at its left  and its right, respectively, that separate it from the closest larger gaps;

\unitlength1.30mm \linethickness{1pt}
\begin{picture}(111.33,12.27)
\put(13.00,5.00){\makebox(0,0)[lc]{$K:$}}
\put(20.00,5.00){\line(1,0){5.00}}
\put(35.00,5.00){\line(1,0){10.00}}
\put(53.00,5.00){\line(1,0){6.00}}
\put(70.00,5.00){\line(1,0){5.00}}

\put(20.00,3.70){\makebox(0,0)[cb]{$|$}}
\put(25.30,3.70){\makebox(0,0)[cb]{$($}}
\put(34.70,3.70){\makebox(0,0)[cb]{$)$}}
\put(45.30,3.70){\makebox(0,0)[cb]{$($}}
\put(52.70,3.70){\makebox(0,0)[cb]{$)$}}
\put(59.30,3.70){\makebox(0,0)[cb]{$($}}
\put(69.70,3.70){\makebox(0,0)[cb]{$)$}}
\put(75.00,3.70){\makebox(0,0)[cb]{$|$}}

\put(40.00,7.00){\makebox(0,0)[cb]{$L_U$}}
\put(49.00,5.50){\makebox(0,0)[cb]{$U$}}
\put(56.00,7.00){\makebox(0,0)[cb]{$R_U$}}
\end{picture}
Put $\tau_R(U):=~\frac{|R_U|}{|U|}$ and $~\tau_L(U):=~\frac{|L_U|}{|U|}$.
The \emph{right thickness}  $\tau _{R}$ and \emph{left thickness} $\tau_{L}$ are
$$\tau_R(K):=~\inf_{U}\tau_R(U),~~~~~~~~~~~\tau_L(K):=~\inf_{U}\tau_L(U)$$
and  Newhouse's thickness $\tau(K)$  is the minimum of $\tau_R(K)$ and $ \tau_L(K)$.

In Section~\ref{sec3},  we show that    for given  two affine  Cantor sets $K$ and $ K'$,  if  $(K,K')$ satisfies the   GTT, then   $K-K'$ is  a   finite union of closed intervals (Proposition \ref{E}).
If  $\tau_R(K)\cdot\tau_L(K')>1$ and  $\tau_L(K)\cdot\tau_R(K')>1$, then   $(K,K')$ satisfies the  SGTT. If  $(K,K')$ satisfies the   SGTT, then   $K$  and $K'$ carry a recurrent set  (Theorem \ref{C}).  The former leads us to introduce the pairs  which have stable intersection, while   do not satisfy the  GTT  (Theorem \ref{D}).

\section{Topological structure of  the sum set  } \label{sec2}
Suppose that  $K$ is a  generalized  homogeneous Cantor set (GHC). Thus,  there are  the rational numbers $\{\frac{r_i}{s_i}\}_{i\in\Sigma}$, $(r_i,s_i)=1$ such that $\log _{p_0}p_i=\frac{\log p_i}{ \log p_0}=\frac{r_i}{s_i}$ for  all  $i\in\Sigma$.  Let $s=\Pi_{  k\in\Sigma}s_k$,   $p:=p_0^{\frac{1}{s}}$ and $m_i:=\frac{r_i}{s_i}s\in\mathbb N$. 
Therefore, we find the real number $p$ and the natural numbers  $\{m_i\}_{i=0}^M$, such that  $p_i=p^{m_i}$.  Moreover,   $\cal{G}$ is dense in  $\cal{A}$, since for each $a>1$,  the set  $\{(a^{r_0},a^{r_1}, \dots ,a^{r_M}) |~ r_i\in\mathbb Q^+, i\in \Sigma\}$ is dense in ${\mathbb R_{>1}}^{M+1}$. Without loss of generality,  we  suppose that $(m_1,m_2,\dots,m_M)=1$ which $(,)$ means   greatest common divisor (gcd)  of the numbers.  
Hence,  affine Cantor set $K$ is a  GHC set, if there are the  real number $p$ and natural numbers $\{m_i\}_{i=0}^M$, such that  $D\psi_{|I_i}=p^{m_i}$.   

\begin{theorem}\label{A}  Suppose that $(K,K')\in  {\cal{A}}\times{\cal{A}}$ and  the set $K-K'$   is not Cantor set. Then  $K-K'$ has dense interior, if at least one of of below conditions holds: 
\begin{itemize}
\item  [$i)$]   $\frac{\log |p_i|}{ \log |p_j'|}\in \mathbb{Q}$ for  all the  pairs $(i,j)\in\Sigma\times\Sigma'$,
\item [$ii)$]   $(K,K')$ belongs to $ {\cal{G}}\times{\cal{G}}$,
\item [$iii)$]  $(K,K')$ belongs to $ {\cal{G}}\times{\cal{A}}$ and $N=1$,
\end{itemize}
\end{theorem}\label{1.1}
\begin{proof}
  Indeed, if $K-K'$ is not a Cantor set, then  $ K- K' $  contains an interval. Thus it is enough  to show that   for each point  $x\in K- K' $ there is  an  interval in $ K- K' $  arbitrarily close to $x$. 
  
 For $(i)$, regard to Theorem  1 in the our paper  \cite{P5}, we introduced an affine map $F:K-K'\to (K-K')\cap (x-R,~x+R)$,  for given $x\in K-K'$ and $0<R<a+b$. The  assertion obtained, since   $(K-K')^\circ\neq \emptyset$.

The assertions $(ii)$ and $(iii)$ are obtained via the following three lemmas.
  
 \begin{lemma}\label{A1}   For each $(K,K')\in  {\cal{A}}\times{\cal{A}}$, if $ K-\lambda K' $  contains an interval, then $ K-\frac{1}{p_i}\lambda K' $ and  $ K-q_j\lambda K' $ contain an interval for all the  pairs $(i,j)\in\Sigma\times\Sigma'$.
\end{lemma}
\begin{proof} Suppose that $I$ is an  interval  in $  K-\lambda K' $ and  $i$ and $j$ are given.  Then, the points $\{\lambda\}\times I$ are difference  pairs.   Consequently,
  all   points of the segments  $X:=T_i^{-1}(\lambda\times I)$  and $Y:={T'_j}^{-1}(\lambda\times I)$   are difference  pairs. This means that     $\pi_2(X)$ and  $\pi_2(Y) $ are  the intervals in  $ K-\frac {1}{p_i}\lambda K' $ and  $K-q_j\lambda K' $ respectively,  which   completes the proof of  the lemma.   For each point  $(x_1,x_2)\in\mathbb{R}^2$  and $k\in \{0,1\}$,  we always get  $ \pi_k(x_1,x_2):=x_k$.
\end{proof} 
 \begin{lemma}\label{A2}   For each $(K,K')\in  {\cal{G}}\times{\cal{G}}$, if $ K-\lambda K' $  contains an interval, then  $ K-p_i\lambda K' $  and $ K-\frac{1}{q_j}\lambda K' $ contain an interval for all the  pairs $(i,j)\in\Sigma\times\Sigma'$.
\end{lemma}
\begin{proof}  Since 
$$ K-\lambda K' =\bigcup_{i\in\Sigma}\pi_2\Big(T_i^{-1}\big(\{p_i\}\times ( K-p_i\lambda K')\big)\Big),$$
from  Baire's theorem, there is an   $i_1\in\Sigma$, such that $ K-p_{i_1}\lambda K' $ contains an interval. Continue this process until you find one sequence $\{i_n\}$ such that $ K-\Pi_{k=1}^np_{i_k}\lambda K' $ contains an interval. 
There exists an element such as $k \in \Sigma$, that appears  in  the sequence $\{i_n\}$ infinitely often.  From Lemma \ref{A1},  $ K-p_k^n\lambda K' $ contains an interval, for all $n\in\mathbb N$. 
There exist a finite sequence $\{\overline{m_i}\}_{k\neq i\in\Sigma}$ and a number $\overline{n}$ of the set ${\mathbb N}_*:=\mathbb N\cup\{0\}$ that satisfy the Diophantine equation:
$$\overline{n}m_k-\sum_{k\neq i\in\Sigma}\overline{m_i}m_i=1.$$
Thus, $ K-p\lambda K' $ contains an interval, since
$$p=\frac{p^{\overline{n}m_k}}{\Pi_{{k\neq i\in\Sigma}}p^{\overline{m_i}m_i}}=\frac{p_k^{\overline{n}}}{\Pi_{{k\neq i\in\Sigma}}p_{m_i}^{\overline{m_i}}}.$$  By induction method, $ K-p^n\lambda K' $ contain an interval, for each $n\in \mathbb N$,  Thus,  $ K-p_i\lambda K' $ contains an interval for all the  $i\in\Sigma$. The second assertion  has a similar proof, and thus the proof of the lemma is complete.
\end{proof} 
\begin{lemma}\label{A3}   If $(K,K')\in  {\cal{A}}\times{\cal{A}}$ and   $N=1$,  then $ K-\frac{1}{q_j}\lambda K' $, $j\in \{0,1\}$ contains an interval.
\end{lemma}
\begin{proof}  Since $$ K-\lambda K' =(K-\frac{\lambda}{q_0} K')\cup\big(K-\frac{\lambda}{q_1}K'-\lambda+\frac{\lambda}{q_1}\big),$$
 one of the sets  $ K-\frac{1}{q_0}\lambda K' $ or  $ K-\frac{1}{q_1}\lambda K' $  contains an interval. We prove that if  $ K-\frac{1}{q_0}\lambda K' $  contains an interval, then   $ K-\frac{1}{q_1}\lambda K' $ also  does. The other case takes the  similar proof.  
Reasoning by contradiction, suppose that $(  K-\frac{1}{q_1}\lambda K' )^\circ=\emptyset$.
 
\textbf{Claim 1}:   For given   numbers   $k\in\mathbb{N}^*$ and  $l\in\mathbb{N}$,  we have $( K-\frac{1}{q_0^kq_1^l}\lambda K' )^\circ=\emptyset$.
For the proof of the claim,   assume that $( K-\frac{1}{q_0^kq_1^l}\lambda K' )^\circ\neq\emptyset$. By using  Lemma \ref{B} $k+l-1$ times, we obtain $(  K-\frac{1}{q_1}\lambda K' )^\circ\neq\emptyset$ which  makes a contradiction and completes the claim. 

For each $j\in\mathbb{N}\cup\{0\}$,  put $t_j:=\frac{b}{q_0^j}(1-\frac{1}{q_1})$. Thus,
$$K'=\{0\}\cup\cup_{j=0}^\infty(\frac{1}{q_0^jq_1}K+t_j).$$ Consequently
$$  K-\frac{1}{q_0}\lambda K'=K\bigcup\cup_{j\in\mathbb N}\big (K-\frac{1}{q_0^jq_1}\lambda K'-\lambda t_j\big).$$ By Claim 1 and Baire's theorem, it follows that $K-\frac{1}{q_0}\lambda K' $ has empty interior.
\end{proof} 
 To proof of  $(ii)$,  suppose that $\epsilon>0.$  Since $(1,x)$ is a difference pair,  there is the finite  sequences $\{i_k\}_{k=1}^{k=m}$  of the set $\Sigma$ and $\{j_k\}_{k=1}^{k=m'}$ of the set $\Sigma'$,  and also natural numbers $R:=\Pi_{k=1}^{k=m} p_{ik}$ and $S:=\Pi_{k=1}^{k=m'} q_{jk}$ such that 
$$
T '_{j_{m'}}\circ\ldots\circ T'_{j_{1}}\circ T
_{i_{m}}\circ\ldots\circ T _{i_{1}}
(1,x)\in \big\{\frac{p^R}{q^S}\big\}\times\big[-\frac{p^Rb}{q^S},~a\big].$$
We also can select the enough big  natural numbers $R,S$ such  that
$$p^R \epsilon>2(a+pqb)~~~~~~~\text{and}~~~~~~~ \frac{p^R}{q^S}\in [1,~pq].$$
On the other hand, based on lemma \ref{A2}, there exists an interval  $I\subset  K- \frac{p^R}{q^S} K' $. Thus,  all points of the segment
$$
J:=\big(T '_{j_{m'}}\circ\ldots\circ T'_{j_{1}}\circ T
_{i_{m}}\circ\ldots\circ T _{i_{1}}\big)^{-1}
\big( \big\{\frac{p^R}{q^S}\big\}\times I\big)
$$ are difference pairs.  Thus, the interval $\pi_2(J)$  has   the sise $\frac{|I|}{p^R}\leq\frac{ a+pqb}{p^R}<\frac{\epsilon}{2}$ and so is subset of $[x-\epsilon,~x+\epsilon]$. This  completes the proof of $(ii)$,  since $\pi_2(J)\subset K-K'$.

For $(iii)$,   suppose that $\epsilon>0.$  Since $(1,x)$ is a difference pair,  there is the finite  sequences $\{i_k\}_{k=1}^{k=m}$  of the set $\Sigma$ and $\{j_k\}_{k=1}^{k=m'}$ of the set $\Sigma'=\{0, 1\}$,  and also natural numbers $R,S$ and $N:=\Pi_{k=1}^{k=m} p_{ik}$ such that 
$$
T '_{j_{m'}}\circ\ldots\circ T'_{j_{1}}\circ T
_{i_{m}}\circ\ldots\circ T _{i_{1}}
(1,x)\in \big\{\frac{p^N}{q_0^Rq_1^S}\big\}\times\big[-\frac{p^Nb}{q_0^Rq_1^S},~a\big].$$
Using an argument similar to the previous case and invoking lemmas \ref{A2} and \ref{A3}, the existence of an interval of $K-K'$ in  $[x-\epsilon, x+\epsilon]$ is established,  and  this  completes the proof of the theorem.
\end{proof}
The following two results can be easily derived from Lemmas \ref{A1}, \ref{A2}, \ref{A3} and their proofs.
\begin{corollary} Suppose that $(K,K')\in {\cal{G}}\times{\cal{G}}$ and  $(m_1,m_2,\dots,m_M)=(n_1,n_2,\dots,n_N)=1$. If $ K- K' $  contains at least one interval, then for all $n,m\in \mathbb{Z}$, the set  $ K-\frac{p^n}{q^m} K' $  has dense interior. 
\end{corollary}

\begin{corollary} Suppose that $(K,K')\in {\cal{G}}\times{\cal{A}}$ with   $(m_1,m_2,\dots,m_M)=1$ and $N=1$. If $ K- K' $  contains at least one interval, then for all $n,m,k\in \mathbb{Z}$, the set  $ K-\frac{p^n}{q_0^mq_1^k} K' $  has dense interior. 
\end{corollary}
Note that for every pair $(K, K') \in \mathcal{A} \times \mathcal{A}$ satisfying the condition $\frac{\log p_i}{\log p_j'} \in \mathbb{Q}$,   for  all the  pairs $(i,j)\in\Sigma\times\Sigma'$,  we have $\frac{\log p_i}{ \log p_j}\in \mathbb{Q}$ for  all the  pairs $(i,j)\in\Sigma\times\Sigma$, and so  $(K,K')\in  {\cal{G}}\times{\cal{G}}$.  But if  $\frac{\log p_i}{ \log p_j}\in \mathbb{Q}$ for  all the  pairs $(i,j)\in\Sigma\times\Sigma$, it is possible that  $\frac{\log p_i}{ \log p_j'}\notin \mathbb{Q}$ for  all the  pairs $(i,j)\in\Sigma\times\Sigma'$, 
\begin{theorem}\label{B}  There is a subset $\cal{F}\subset {\cal{N}}\times{\cal{N}}$ such that for  each  $(K,K')\in\cal{F}$, $K-K'$ is  a bilateral gap interval  an  L, R, M-Cantorval or Cantor set. The set $\cal{F}$ has full size and contains  an open  and dense subset. 
\end{theorem}\label{1.1}
\begin{proof} The assertion  is obtained  via the following three lemmas.
\begin{lemma}\label{B1} Suppose that $1<\gamma$ be a real number such that  $p_M=\gamma^{m}$ and $q_0=\gamma^{n}$ with $(m,n)=1$.  Suppose that there are $\lambda, x$ and $\epsilon>0$ with  
$$\max\{a+\epsilon+\frac{b-bq_N}{q_N}(\lambda-\epsilon), a-\frac{a}{p_M}\}<x<a$$ such for all $(s,t)\in \mathbb B_{\epsilon}(\lambda)\times\mathbb B_{\epsilon}(x)$, the pair $(s,t)$ is  not a difference pair. Thus,  if $\gamma<\frac{\lambda+\epsilon}{\lambda-\epsilon}$, then   for each  $\mu>0$ the point   $a$ is accumulated  by infinitely many  gaps of $K-\mu K'$.
\end{lemma}
\begin{proof}   None of the  points  of the regions  $U_k:={T_M}^{-nk}\circ {T_0'}^{-mk}\big(B_{\epsilon}(\lambda)\times\mathbb B_{\epsilon}(x)\big),~k\in\mathbb N$ (which its boundary is indeed a rectangle) are difference pairs, since $U_1$ and so the all $U_k$s always  stay above the lines $L_1^N:~t=\frac{b-bq_N}{q_N}s+a$ and $t=a-\frac{a}{p_M}.$   Thus,   for each  $\mu\in B_{\epsilon}(\lambda)$ the point   $a$ is accumulated  by infinitely many  gaps of $K-\mu K'$, since
${T_M}^{-nk}\circ {T_0'}^{-mk}(s,t)=\big(s,p_M^{-nk}t+a(1-p_M^{-nk})\big)$ for all  $k\in\mathbb N$. 
 For all $i,j\in\mathbb N$ satisfying $1\leq \frac{q_0^{j}}{p_M^{i}} \leq p_Mq_0,$ let the intervals
$$I_{i,j}:=\pi_1\circ{T_M}^{-i}\circ {T_0'}^{-j}\big(B_{\epsilon}(\lambda)\times\mathbb B_{\epsilon}(x)\big),$$  
where   $\pi_1$ is the map  that  projects the points of $\mathbb{R}^2$   to their first component.
Obviously, $I_{i+n,j+m}=I_{i,j}$ for each $i$ and  $j$. Also for each $k\in \mathbb N$, Diophantine equation $nj-mi=k$ always  has the solution, since $(m,n)=1$.  Thus,
 $$\bigcup_{nj-mi=0} ^{nj-mi=m+n}I_{i,j}= \bigcup_{nj-mi=0} ^{nj-mi=m+n}\frac{q_0^{j}}{p_M^{i}}B_{\epsilon}(\lambda)=\bigcup_{k=0} ^{m+n}\gamma^{k}B_{\epsilon}(\lambda).$$
Above set is connected if and only if $\gamma(\lambda-\epsilon)<\lambda+\epsilon$. Hence, $\cup I_{i,j}=[\lambda-\epsilon,~p_Mq_0(\lambda+\epsilon)]$, since  $\gamma<\frac{\lambda+\epsilon}{\lambda-\epsilon}$.  Consequently, non of  the elements  of $U_{i,j}= {T_M}^{-i}\circ{T_0'}^{-j}\big(B_{\epsilon}(\lambda)\times\mathbb B_{\epsilon}(x)\big)$, where $0\leq nj-mi \leq m+n$ are  difference pairs. From  the above proof, for each  $\lambda-\epsilon <\mu<p_Mq_0(\lambda+\epsilon)$ and so  each  $\mu>0$   the point   $a$ is accumulated  by infinitely many  gaps of $K-\mu K'$ and this completes the proof.
\end{proof}\label{E} 
From  proof of Lemma 2,  one can prove that  if $\frac{\log p_1}{ \log q_0}\notin \mathbb{Q}$ and    the pair $(s,t)$ with $\max\{a+\frac{b-bq_N}{q_N}s, a-\frac{a}{p_1}\}<t<a$  is  not a difference pair, then   for each  $\mu>0$ the point   $a$ is accumulated  by infinitely many  gaps of $K-\mu K'$.
Let,
$${\cal{B}}:=\big\{(K,K')\in {\cal{N}}\times {\cal{N}}\big| ~\forall \lambda>0, ~\exists\delta>0,~s.t~ [a-\delta,~a]\subset K-\lambda K'\big\}.
$$
\begin{lemma}\label{B2} Suppose that $(K,K')\in {\cal{B}}$ and  $(c,d)$ be a bounded gap of $K-K'$. Then
 there is  a  $\delta>0$ such  that $[c-\delta,~c]\subset K- K'$.  

\end{lemma}\label{G}
\begin{proof} 
There are  the  finite sequences  $\{a_{i}\}^{m}_{i=1},~\{b_{j}\}^{n}_{j=1}$ of the set $\Sigma$ and $\Sigma'$, respectively,  and  also $\lambda>0$ such that 
$$T' _{b_{n}}\circ T'_{b_{n-1}} \circ\ldots\circ
T'_{b_{1}}\circ T_{a_{m}}\circ T _{a_{m-1}} \circ\ldots\circ
T_{a_{1}}(1,c)=(\lambda,a).$$
Thus,  there is a $\beta>0$ such that $[c-\beta,~c]\subset K-K'$ since  there is a $\delta>0$ such that $ [a-\delta,~a]\subset K-\lambda K'$.
\end {proof} 
\begin{lemma}\label{B3}
 There is an open and dense full measure  subset ${\cal{U}}$  of $ {\cal{N}}\times {\cal{N}}-{\cal{B}}$, such that if 
 $(K,K')\in {\cal{U}} $ and  $(c,d)$ be a bounded gap of $K-K'$, then  the point   $c$ is accumulated  by infinitely many  gaps of $K-K'$.
\end{lemma}\label{G}
\begin{proof} Fix two Cantor sets  $K=K_{p_0,p_1,\dots,p_M}^{a_1,a_2,\dots,a_{M-1},a}\in  {\cal{N}}$ and $K'={K'}_{q_0,q_1,\dots,q_N}^{b_1,b_2,\dots,b_{N-1},b}\in  {\cal{N}}$ with  the pair $(K,K')\notin {\cal{B}}$, then  there is a $\lambda>0$  such that $K-\lambda K'$   has a gap near to $a$.  For  $\epsilon>0$ small, let $V_{\epsilon}^{(K,K')}$  be the set of all $(\tilde K,\tilde {K'})\in  {\cal{N}}\times {\cal{N}}$ such that if  $\tilde{K}=\tilde{K}_{\tilde{p_0},\tilde{p_1},\dots,\tilde{p_M}}^{\tilde{a_1},\tilde{a_2},\dots,\tilde{a_{M-1}},\tilde{a}}$ and $\tilde{K'}={\tilde{K'}}_{\tilde{q_0},\tilde{q_1},\dots,\tilde{q_N}}^{\tilde{b_1},\tilde{b_2},\dots,\tilde{b_{N-1}},\tilde{b}}$, then 
 $$(\tilde{p_0},\tilde{p_1},\dots,\tilde{p_M},\tilde{q_0},\tilde{q_1},\dots,\tilde{q_N},\tilde{a},\tilde{\lambda})\in  \Pi_{i=0}^{M} \mathbb{B}_\epsilon(p_i)\times \Pi_{j=0}^{N}\mathbb{B}_\epsilon(q_j)\times \mathbb{B}_\epsilon (a)\times \mathbb{B}_\epsilon(\lambda),$$ and 
 $$(\tilde{a_1},\tilde{a_2},\dots,\tilde{a_{M-1}},\tilde{b_1},\tilde{b_2},\dots,\tilde{b_{N-1}})\in  \Pi_{i=1}^{M-1} \mathbb{B}_\epsilon(a_i)\times \Pi_{j=1}^{N-1}\mathbb{B}_\epsilon(b_j).$$ 
 One can find $\epsilon>0$ small and a $x$ near to $a$, such that for each  $(\tilde K,\tilde {K'})\in V_{\epsilon}^{(K,K')}$ the set  $K_{\tilde{p_0},\tilde{p_1},\dots,\tilde{p_M}}^{\tilde{a_1},\tilde{a_2},\dots,\tilde{a_{M-1}}\tilde{a}}-\tilde{\lambda} {K'}_{\tilde{q_0},\tilde{q_1},\dots,\tilde{q_N}}^{\tilde{b_1},\tilde{b_2},\dots,\tilde{b_{N-1}},b}$ has a gap near to $\tilde{a}$ which contains $\mathbb{B}_\epsilon (x)$  (we need that  $\max\{a+2\epsilon+\frac{b-b(q_N-\epsilon)}{q_N-\epsilon}(\lambda-\epsilon), a+\epsilon-\frac{a+\epsilon}{p_M+\epsilon}\}<x<a-\epsilon$).   Let $\gamma_0=\frac{\lambda+\epsilon}{\lambda-\epsilon}$ and  let $U_{\epsilon}^{(K,K')}$  be the set of all $(\tilde K,\tilde {K'})\in V_{\epsilon}^{(K,K')}$ such that
$$(p_M,q_0)\in \mathbb{B}_\epsilon(p_M)\times \mathbb{B}_\epsilon(q_0)-\{(\gamma ^i,\gamma ^j)|\gamma\geq \gamma_0,~i,j\in\mathbb N\}.$$
For each the pair  $(\tilde K,\tilde {K'})\in U_{\epsilon}^{(K,K')}$,
\begin{itemize}
\item  If  $\frac{\log {\tilde{p}_M}}{ \log \tilde{q_0}}\notin \mathbb{Q}$, then for each  $\mu>0$,  the point  $\tilde{a}$ is accumulated  by infinitely many  gaps of  $K_{\tilde{p_0},\tilde{p_1},\dots,\tilde{p_M}}^{\tilde{a_1},\tilde{a_2},\dots,\tilde{a_{M-1}}\tilde{a}}-\tilde{\mu} {K'}_{\tilde{q_0},\tilde{q_1},\dots,\tilde{q_N}}^{\tilde{b_1},\tilde{b_2},\dots,\tilde{b_{N-1}},b}$.
\item  If  $\frac{\log {\tilde{p}_M}}{ \log \tilde{q_0}}\in \mathbb{Q}$, then there is a   $\gamma<\gamma_0$ such that  $\tilde{p_M}=\gamma^{n}$ and $\tilde{q_0}=\gamma^{m}$ where $(m,n)=1$.  Since  $\max\{\tilde{a}+\epsilon+\frac{b-b\tilde{q_N}}{\tilde{q_N}}(\lambda-\epsilon), \tilde{a}-\frac{\tilde{a}}{\tilde{p}_M}\}<x<\tilde a$ and  for all $(s,t)\in \mathbb B_{\epsilon}(\lambda)\times\mathbb B_{\epsilon}(x)$, the pair $(s,t)$ is  not a difference pair for $K_{\tilde{p_0},\tilde{p_1},\dots,\tilde{p_M}}^{\tilde{a_1},\tilde{a_2},\dots,\tilde{a_{M-1}}\tilde{a}}$  and $  {K'}_{\tilde{q_0},\tilde{q_1},\dots,\tilde{q_N}}^{\tilde{b_1},\tilde{b_2},\dots,\tilde{b_{N-1}},b}$ . From Lemma 2,  for each  $\mu >0$,  the point  $\tilde{a}$ is accumulated  by infinitely many  gaps of  $K_{\tilde{p_0},\tilde{p_1},\dots,\tilde{p_M}}^{\tilde{a_1},\tilde{a_2},\dots,\tilde{a_{M-1}}\tilde{a}}-\tilde{\mu} {K'}_{\tilde{q_0},\tilde{q_1},\dots,\tilde{q_N}}^{\tilde{b_1},\tilde{b_2},\dots,\tilde{b_{N-1}},b}$.
\end{itemize}

Obviously, $U_{\epsilon}^{(K,K')}$ is an open and dense full measure subset of $V_{\epsilon}^{(K,K')}$, since $\gamma_0>1$.   Note that it is possible that $(K,K')\notin U_{\epsilon}^{(K,K')}$. Let,
$${\cal{U}}=\bigcup _{(K,K')\in { {\cal{N}}\times {\cal{N}}}-{\cal{B}}}U_{\epsilon}^{(K,K')}$$
which is  so  an open and dense full measure  subset of $ {\cal{N}}\times {\cal{N}}-{\cal{B}}$. 
Also for all  $(K,K')\in {\cal{U}}$ and  $\mu>0$, the point   $a$ is accumulated  by infinitely many  gaps of $K-\mu K'$.

Suppose that  $(K,K')\in {\cal{U}}$ and $\{(1,x_i)\}$ is   a sequence   of  the non difference  pairs  such that $x_i\to a$. There is a $\epsilon>0$ such that for each $i\in\mathbb N$, the pairs $(s,t)\in\mathbb{B}_{\epsilon} (1)\times \mathbb{B}_{\epsilon_i}(x_i)$ are   not  difference pair.  Thus, $\epsilon_i\to 0$ when $i\to\infty$.
For all  the finite  sequences $\{a_i\}_{i=1}^{i=m}$ and $\{b_j\}_{j=1}^{j=n}$ of the set $\Sigma$ and $\Sigma'$,  respectively,  denote
$$T_{a_1,\dots,a_m}^{b_1,\dots,b_{n}}:=T '_{b_{n}}\circ\ldots\circ T'_{b_{1}}\circ T
_{a_{m}}\circ\ldots\circ T _{a_{1}}.$$

There are the numbers  $\lambda_{a_1,\dots,a_m}^{b_1,\dots,b_{n}}:=\frac{\Pi_{k=1}^{m}p_{a_k}}{\Pi_{k=1}^{n}q_{b_k}}$ and  $t_{a_1,\dots,a_m}^{b_1,\dots,b_{n}}$ such that 
 $$T_{a_1,\dots,a_m}^{b_1,\dots,b_{n}}(s,t)=(\lambda_{a_1,\dots,a_m}^{b_1,\dots,b_{n}} s, p_{a_1}\dots p_{a_m}t+ t_{a_1,\dots,a_m}^{b_1,\dots,b_{n}}).$$
Thus, for each $x\in K-K'$, we always have 
$$T_{a_1,\dots,a_m}^{b_1,\dots,b_{n}}
(1,x)\in \big\{\lambda_{a_1,\dots,a_m}^{b_1,\dots,b_{n}}\big\}\times\big[-\lambda_{a_1,\dots,a_m}^{b_1,\dots,b_{n}}b,~a\big].$$
Fix $k\in\mathbb N$ and let
$$\tilde{E_k}:=\{T_{c_1,\dots,c_m}^{d_1,\dots,d_{n}}|~~~1\leq\lambda_{c_1,\dots,c_m}^{d_1,\dots,d_{n}}< pq,~~ p_{a_1}\dots p_{a_m}>k\}.$$
 Remove elements of $\tilde{E_k}$ as follows: if $T\neq T'\in \tilde{E_k}$  and  there are  finite  sequences  $\{a_i\}_{i=1}^{i=m}$ and $\{b_j\}_{j=1}^{j=n}$  of the set $\Sigma$ and $\Sigma'$, respectively such that $T'= T_{a_1,\dots,a_m}^{b_1,\dots,b_{n}}\circ T$, then remove $T'$. Denote  such set to $E_k$, which is obviously finite and  
$$K- K'=\pi_2\Big( \bigcup_{T\in E_k} T^{-1}\big(\big\{\lambda_{a_1,\dots,a_m}^{b_1,\dots,b_{n}}\big\}\times \big(K-\lambda_{a_1,\dots,a_m}^{b_1,\dots,b_{n}}K')\big)\Big).$$
where   $\pi_2$ is the map  that  projects the points of $\mathbb{R}^2$   to their second component.
Let $J_{a_1,\dots,a_m}^{b_1,\dots,b_{n}}=\pi_2\Big( T^{-1}\big(\big\{\lambda_{a_1,\dots,a_m}^{b_1,\dots,b_{n}}\big\}\times \big[-\lambda_{a_1,\dots,a_m}^{b_1,\dots,b_{n}}b,~a\big]\big)\Big)$ and 
${\cal{I}}_k=\{J_{a_1,\dots,a_m}^{b_1,\dots,b_{n}}|~T_{a_1,\dots,a_m}^{b_1,\dots,b_{n}}\in E_k\}$.
Since $(c,d)$ is a bounded gap of $K-K'$ and $p_{a_1}\dots p_{a_m}>k$ for all $T\in E_k$, there is a $k\in\mathbb N$ such that $d-c>\frac{a+pqb}{p_{a_1}\dots p_{a_m}}$   for all $T\in E_k$. Note that $c$ may be equal to $\max I$ for more than one of the intervals $I \in \mathcal{I}_k$.
Let $A_k=\{J\in {\cal{I}}_k| ~c=\max I\}$ and
$$\alpha=\min_{I\notin A_k}\{c-  \max I|~~  \max I<c\}.$$
Consequently,
$$K- K'=\pi_2\Big( \bigcup_{I_{a_1,\dots,a_m}^{b_1,\dots,b_{n}}\in A_k} T^{-1}\big(\big\{\lambda_{a_1,\dots,a_m}^{b_1,\dots,b_{n}}\big\}\times \big(K-\lambda_{a_1,\dots,a_m}^{b_1,\dots,b_{n}}K')\big)\Big),$$
 on the  $(c-\alpha,~c]$. Since $\pi_2\circ T_{a_1,\dots,a_m}^{b_1,\dots,b_{n}}
(1,c)=a$, where $J_{a_1,\dots,a_m}^{b_1,\dots,b_{n}}\in A_k$, one may choose $0<\beta<\alpha$, such that 
 $ T_{a_1,\dots,a_m}^{b_1,\dots,b_{n}}
(1,t)$  always  stay above the lines $L_1^N:~t=\frac{b-bq_N}{q_N}s+a$ and $t=a-\frac{a}{p_M}$,   for all $t\in (c-\beta,~c]$. This indicates that only the operators $T_M$ and $T'_0$ can be employed to make $T_{a_1,\dots,a_m}^{b_1,\dots,b_{n}}(1,t)$ a difference pair.
 For each $J_{a_1,\dots,a_m}^{b_1,\dots,b_{n}}\in A_k$, let $m_{a_1,\dots,a_m}^{b_1,\dots,b_{n}}\in \mathbb N$
and  $n_{a_1,\dots,a_m}^{b_1,\dots,b_{n}}\in \mathbb N$  such that $p^{m_{a_1,\dots,a_m}^{b_1,\dots,b_{n}}}=\Pi_{k=1}^{m}p_{a_k}$ and $q^{n_{a_1,\dots,a_m}^{b_1,\dots,b_{n}}}=\Pi_{k=1}^{n}q_{b_k}$. 
 Order the set $\{\frac{p^r}{q^s}|~ 1\leq\frac{p^r}{q^s}\leq pq,~~~~r,s\in\mathbb N\}$ as follows: $\frac{p^{r_1}}{q^{s_1}}<\frac{p^{r_2}}{q^{s_2}}$, if $r_1<r_2$, $r_1=r_2$ then $s_2<s_1$ (do not simplify the fractions). If  $\frac{p^{r_1}}{q^{s_1}}<\frac{p^{r_2}}{q^{s_2}}$, denote $T_{\frac{p^{r_1}}{q^{s_1}},\frac{p^{r_2}}{q^{s_2}}}:={T'_0}^{s_2-s_1}\circ T_M^{r_2-r_1}$. Take $m,n\in \mathbb N$ such that  $\frac{p^m}{q^n}<\epsilon$ and with  above order 
 $$\frac{p^m}{q^n}>\max\{\frac{p^{m_{a_1,\dots,a_m}^{b_1,\dots,b_{n}}}}{q^{n_{a_1,\dots,a_m}^{b_1,\dots,b_{n}}}}|~~J_{a_1,\dots,a_m}^{b_1,\dots,b_{n}}\in A_k\}.$$
 
 The operators  $S_{a_1,\dots,a_m}^{b_1,\dots,b_{n}}$, where $J_{a_1,\dots,a_m}^{b_1,\dots,b_{n}}\in A_k$  given by
 $$S_{a_1,\dots,a_m}^{b_1,\dots,b_{n}}(s,t):=T_{\frac{p^m}{q^n},\frac{p^{m_{a_1,\dots,a_m}^{b_1,\dots,b_{n}}}}{q^{n_{a_1,\dots,a_m}^{b_1,\dots,b_{n}}}}}\circ T_{a_1,\dots,a_m}^{b_1,\dots,b_{n}}
(s,t):=(\frac{p^m}{q^n}s,p^mt+z_{a_1,\dots,a_m}^{b_1,\dots,b_{n}})$$ are an operator, since $S_{a_1,\dots,a_m}^{b_1,\dots,b_{n}}(1,c)=(\frac{p^m}{q^n}, a)$. Let $t_i\in (c-\beta,~c]$ such that  $S_{a_1,\dots,a_m}^{b_1,\dots,b_{n}}(1,t_i)=(\frac{p^m}{q^n}, x_i)$, for one of  $J_{a_1,\dots,a_m}^{b_1,\dots,b_{n}}\in A_k$.
 Thus, the pairs  $\{(1,x_j)|~j\geq i\}$are  not the difference pairs for $K$  and $K'$. This tells that $x_j\notin K-K'$,  for all $j\geq i$, and so  the point   $c$ is accumulated  by infinitely many  gaps of $K- K'$.
\end{proof} 
 Let ${\cal{C}}:={\cal{U}}\bigcup {{\cal{B}}}$, where ${\cal{U}}$ is the set  appeared in Lemma \ref{B3}.  
It is easy to check that ${\cal{C}}$  has full measure in  $\cal{N}\times \cal{N}$ and   ${{\cal{C}}}^\circ$ is  a  dense  subset  of  ${\cal{C}}$.   To proof of density, from above proof, $ {\cal{B}}$ is closed in  $\cal{N}\times \cal{N}$ and so  ${\cal{U}}$ is open in $\cal{N}\times \cal{N}$. This says  ${\cal{U}}\bigcup {{\cal{B}}}^\circ={({\cal{U}}\bigcup {\cal{B}})}^\circ$. Now if  $\cal W$ is an open set containing $(K,K')\in {\cal{B}}$,  two cases take place.  If   ${\cal W}\subset {\cal B}$, then  $(K,K')\in {\cal{B}}^\circ$.  Otherwise, there is a   $(\tilde K,\tilde {K'})\in  {\cal{W}}- {\cal{B}}$ and  $\epsilon>0$ small  such that  $U_{\epsilon}^{(\tilde K,\tilde {K'})}\subset {\cal{W}}- {\cal{B}}$. Thus,  ${\cal W}\cap {\cal U}\neq\emptyset$, and this completes   the claim.

For each  $X\subset {\cal{N}}\times {\cal{N}}$, denote $\tilde{X}:=\{(K,K')|~(K',K)\in X\}$.
Now, we are ready to  prove Theorem 2.  

Let $\cal{F}=\cal{C}\cap\tilde{\cal{C}}$  which  obviously    has full measure  and contains an open and dense  subset of  $\cal{N}\times \cal{N}$. 
 Suppose that  $(K,K')\in\cal{F}$, and $K-K'$ is not Cantor set. Thus  for each point  $x\in K- K' $, there is  an  interval in $ K- K' $  arbitrarily close to $x$, Hence,

\begin{itemize}
\item  If $(K,K')\in {\cal{U}}\cap\tilde{\cal{U}}$, then  from Lemma 4, $K-K'$ is an  M-Cantorval,
\item  If $(K,K')\in {\cal{U}}\cap\tilde{\cal{A}}$,  then  from Lemmas 3 and  4, $K-K'$ is an  L-Cantorval,
\item  if $(K,K')\in {\cal{A}}\cap\tilde{\cal{U}}$,  then   from Lemmas 3 and  4, $K-K'$ is an R-Cantorval,
\item  if $(K,K')\in {\cal{A}}\cap\tilde{\cal{A}}$, then  each gap of  $K-K'$ located between too intervals,  if it is not an interval.
\end{itemize}
 These facts completes the proof of Theorem \ref{B}.
 \end{proof} 
The following theorem states that there are no examples of stable, orientation-preserving expansive mappings that lack the five aforementioned properties \cite {M4}.
 \begin{theorem}\label{BB}  There is a subset ${\cal{D}}\subset \{(K,K')\in {\cal{A}}\times{\cal{A}}|~ D\psi,D\psi'>0\}$ such that for  each  $(K,K')\in\cal{D}$, $K-K'$ is  a bilateral gap interval  an  L, R, M-Cantorval or Cantor set. 
\end{theorem}
\begin{proof}  Let 
 $${\cal{D}}:=\big\{(K,K')\in{\cal{A}}\times{\cal{A}}\big| ~p=q, ~~~(m_M,n_0)=1,~(m_0,n_N)=1\big\},$$
  Obviosly  $\cal{D}$ is dense in $ \{(K,K')\in {\cal{A}}\times{\cal{A}}|~ D\psi,D\psi'>0\}$.   Suppose that $(K,K')\in {\cal D}$ and let
  $$A:=\{p^{rm_M-sn_0}|~r,s\in\mathbb N_*,~~~0\leq rm_M-sn_0\leq  m_M+n_0\}=\{1,p,p^2,\dots,p_Mq_0\}.$$ If $(1,x_1)$ is not  a difference pair for $K$ and $K'$,  then $\{(1,x_i)\}$ where $x_i=\frac{x_1+p^{in_0m_M}a-a}{p^{in_0m_M}}$,  is   a sequence   of  the non difference  pairs  and $x_i\to a$, Note that,   $p_M^{in_0}=q_0^{im_M}$ and  $T_M^m(s,t)=(p_M^ms,p_M^mt-p_M^ma+a), ~{T'}_0^n(s,t)=(\frac{s}{q_0^n},t)$, for all $m,n\in\mathbb N$. This says $\frac{x_i-a}{x_j-a}=p^{(j-i)n_0m_M}$, for all $i,j\in\mathbb N$. In the other word, for each $i,k\in\mathbb N$, there is a $j\in\mathbb N$ such that  $x_j-a=\frac{x_i-a}{p^{kn_0m_M}}$. 
  Let $r_0,s_0\in\mathbb N_*$ such that $r_0m_M-s_0n_0=-1$, and then  let
 \begin{eqnarray*}
T_0(s,t)={T'_0}^{s_0}\circ {T_M}^{r_0}(s,t)&=&(\frac{p^{r_0m_M}}{p^{s_0n_0}}s, p^{r_0m_M}t- p^{r_0m_M} a_M )\\
&=&(\frac{s}{p}, p^{r_0m_M}t- p^{r_0m_M} a+a).
\end{eqnarray*}
 Considering the notations of lemma \ref{B3},  we obtain $\frac{p^{m_{a_1,\dots,a_m}^{b_1,\dots,b_{n}}}}{q^{n_{a_1,\dots,a_m}^{b_1,\dots,b_{n}}}}=p^{m_{a_1,\dots,a_m}^{b_1,\dots,b_{n}}-n_{a_1,\dots,a_m}^{b_1,\dots,b_{n}}}\in A$,   for all  $J_{a_1,\dots,a_m}^{b_1,\dots,b_{n}}\in A_k$.  For the operators  $S_{a_1,\dots,a_m}^{b_1,\dots,b_{n}}$, where $J_{a_1,\dots,a_m}^{b_1,\dots,b_{n}}\in A_k$  given by
 $$S_{a_1,\dots,a_m}^{b_1,\dots,b_{n}}:={T_0}^{m_{a_1,\dots,a_m}^{b_1,\dots,b_{n}}-n_{a_1,\dots,a_m}^{b_1,\dots,b_{n}}}\circ T_{a_1,\dots,a_m}^{b_1,\dots,b_{n}}
,$$ there are constants $k_{a_1,\dots,a_m}^{b_1,\dots,b_{n}}\in\mathbb N$  and $z_{a_1,\dots,a_m}^{b_1,\dots,b_{n}}$ such that
$$S_{a_1,\dots,a_m}^{b_1,\dots,b_{n}}(s,t)=(s,p^{m_Mn_0k_{a_1,\dots,a_m}^{b_1,\dots,b_{n}}}t+z_{a_1,\dots,a_m}^{b_1,\dots,b_{n}}).$$
 Since $S_{a_1,\dots,a_m}^{b_1,\dots,b_{n}}(1,c)=(1, a)$, we have $z_{a_1,\dots,a_m}^{b_1,\dots,b_{n}}=a-p^{m_Mn_0k_{a_1,\dots,a_m}^{b_1,\dots,b_{n}}}c$. Thus,
 $$R_{a_1,\dots,a_m}^{b_1,\dots,b_{n}}(t):=\pi_2\circ S_{a_1,\dots,a_m}^{b_1,\dots,b_{n}}(s,t)=p^{m_Mn_0k_{a_1,\dots,a_m}^{b_1,\dots,b_{n}}}(t-c)+a.$$ This says
 $$\frac{R_{a_1,\dots,a_m}^{b_1,\dots,b_{n}}(t)-a}{R_{a'_1,\dots,a'_m}^{b'_1,\dots,b'_{n}}(t)-a}=p^{m_Mn_0(k_{a_1,\dots,a_m}^{b_1,\dots,b_{n}}-k_{a'_1,\dots,a'_m}^{b'_1,\dots,b'_{n}})},$$
 for all  $J_{a_1,\dots,a_m}^{b_1,\dots,b_{n}}$ and $J_{a'_1,\dots,a'_m}^{b'_1,\dots,b'_{n}}$ belong to $ A_k$. 
   Let $t_i\in (c-\beta,~c]$ such that  $S_{a_1,\dots,a_m}^{b_1,\dots,b_{n}}(1,t_i)=(1, x_i)$, for one of  $J_{a_1,\dots,a_m}^{b_1,\dots,b_{n}}\in A_k$ and so, $R_{a_1,\dots,a_m}^{b_1,\dots,b_{n}}(t)=x_i$.
   Now, if  $J_{a'_1,\dots,a'_m}^{b'_1,\dots,b'_{n}}\in A_k$, then there is a $j\in \mathbb N$, such that 
   $$R_{a'_1,\dots,a'_m}^{b'_1,\dots,b'_{n}}(t_i)-a=\frac{x_i-a}{p^{m_Mn_0(k_{a_1,\dots,a_m}^{b_1,\dots,b_{n}}-k_{a'_1,\dots,a'_m}^{b'_1,\dots,b'_{n}})}}=x_j-a.$$
   Consequently, $R_{a'_1,\dots,a'_m}^{b'_1,\dots,b'_{n}}(t_i)=x_j$, and we have $S_{a'_1,\dots,a'_m}^{b'_1,\dots,b'_{n}}(1,t_i)=(1,x_j)$. Thus, the pair $(1,t_i)$ is   not a difference pair for $K$  and $K'$. This tells that $t_i\notin K-K'$,  for all $t_i$  enough  near to $c$, and so  the point   $c$ is accumulated  by infinitely many  gaps of $K- K'$. This says   $K- K'$ is an  L-Cantorval.\\
    If  there is $\delta>0,$ such that $ [a-\delta,~a]\subset K-K'$. Fixed a   $J_{a_1,\dots,a_m}^{b_1,\dots,b_{n}}\in A_k$, thus  the all of the points of ${S_{a_1,\dots,a_m}^{b_1,\dots,b_{n}}}^{-1} \big (\{1\}\times [a-\delta,~a]\big)$ are the  difference pairs for $K$  and $K'$. Thus, there is $\delta_1>0,$ such that $ [c-\delta_1,~c]\subset K-K'$, for all gaps $(c,~d)$ of  $K- K'$.
    
    Similar, if   the point   $b$ is accumulated  by infinitely many  gaps of $K- K'$, then  $K- K'$ is an  R-Cantorval. Otherwise,  there is $\delta_2>0,$ such that $ [d, ~d+\delta_2]\subset K-K'$, for all gaps  $(c,~d)$ of  $K- K'$.
    Now, if $K - K'$ is not a Cantor set, then depending on which of the above cases occurs, the set $K - K'$ takes the form of either a  bilateral gap interval or an L, R,  M-type Cantorval, and this  completes the proof of Theorem \ref{BB}.
\end{proof}

\section{Stable intersection of affine Cantor sets} \label{sec3}
The  concepts of  the  lateral thicknesses and   the Generalized Thickness Test  introduced first by   Moreira  as   methods to  check  the  stable intersection of  regular  Cantor sets \cite{M} to   investigated  non hyperbolic diffeomorphisms and  interesting dynamics  created around it.
We begin this section by posing this  two concepts.

Assume  that $K$ and $K'$ are two affine Cantor sets defined by  Markov partitions 
$P:=\cup_{i=0}^{i=M}[a_i,~b_i]$ and $P':=\cup_{j=0}^{j=N}[a'_j,~b'_j]$, respectively. Here,  $I_i, ~i\in \Sigma$  (resp. $I'_j,~j\in\Sigma'$) are not necessarily  of  equal size. 
\begin{definition}
 We say that $(K,K')$ satisfies the  (strictly) Generalized Thickness Test or simply GTT (SGTT), if for each $t\in\mathbb{R}$ and $s>0$ either $P$ is contained in (closure of) one  gap of $s P'+t$  or $s P'+t$  is contained in  (closure of) one  gap of $P$    or $P\cap (s P'+t)\neq\emptyset$ ($P\cap (s P'+t)$ has nonempty interior).
\end{definition}
Two generalization of the Gap Lemma  which has been proved  in \cite{M} are  as below:
\begin{itemize}
\item [$\bullet$] Suppose  that  $\tau_R(K)\cdot\tau_L(K')\geq1$ and  $\tau_L(K)\cdot\tau_R(K')\geq1$, then,  either $K$  is contained in  a gap of $K'$ or $K'$ is contained in a gap of $K$   or $K\cap K'\neq\emptyset$.
\item [$\bullet$] Suppose  that $(K,K')$ satisfies the Generalized Thickness Test, then,  either $K$  is contained in  a gap of $K'$ or $K'$ is contained in a gap of $K$   or $K\cap K'\neq\emptyset$.
\end{itemize}
 
The locus of all points $(s,t),~s>0$ in the plane  $SoT$ that the line $L:~y=\frac{1}{s}x+t, ~s>0$ intercepts $[0,~a]\times [0,~b]$  in the plane $XoY$ is  $ V:=\{(s,t)| ~s>0,~-bs\leq t\leq a\}$. Foe each  ${(i,j)\in\Sigma\times\Sigma'}$, let
$$B_i:= \big\{ (s,t)\big|  ~a_i -bs\leq t\leq b_i\big\}~~\text{and}~~ C_j:=\big\{ (s,t)\big| ~  -b'_js\leq t\leq a -a'_js\big\}.$$
\begin{itemize}
\item If $ \frac{a}{b}\leq s$,  the set $s P'+t$  is  contained in  any  gap of $P$. The set $P$  or  indeed  the interval $[0,~a]$ is not  in a gap $s P'+t$  if and only if $$t\notin [0,~a]-sP'=\cup_{j\in\Sigma'} [-b'_js,~a-a'_js]=\cup_{j\in\Sigma'}C_j.$$
\item If $s<\frac{a}{b} $,  the set $P$ is  contained in  any  gap of $sP'+t$. The set $sP'+t$  or  indeed  the interval $s[0,~b]+t$ is not  in a gap $P$  if and only if $$t\notin P-s [0,~b]=\cup_{i\in\Sigma} [a_i-sb,~b_i]=\cup_{i\in\Sigma}B_i.$$
\end{itemize}

Consequently, $(s,t)\in\bigcup_{{(i,j)\in\Sigma\times\Sigma'}}(B_i\cup C_j)$, if and only if 
$P$ and $s P'+t$ does not  fall in  the each other gaps

\begin{theorem} \label{C}   Assume  that  $K$ and $ K'$  are two affine  Cantor sets. Then 
\begin{itemize}
\item [$i)$]  if  $\tau_R(K)\cdot\tau_L(K')>1$ and  $\tau_L(K)\cdot\tau_R(K')>1$, then   $(K,K')$ satisfies the  SGTT,
\item [$ii)$] 
if  $(K,K')$ satisfies the   SGTT, then   $K$  and $K'$ carry a recurrent set,
\end{itemize}
\end{theorem} \label{AABB} 
\begin {proof}
For $(i)$,  fixed   $s$ and $t$ and suppose  that  $P$ and $s P'+t$ does not  fall in the closure of  the each other gaps. Take $\epsilon>0$ small such that 
$\tau_R(P_1)\cdot\tau_L(P_1')>1$ and  $\tau_L(P_1)\cdot\tau_R(P_1')>1$, where  
$P_1:=\cup_{i=0}^{i=M}[a_i+\epsilon,~b_i-\epsilon]$ and $P_1':=\cup_{j=0}^{j=N}[sa'_j+t+\epsilon,~sb'_j+t-\epsilon]$. Thus $P_1$ and $P'_1$ does not  fall in the closure of  the each other gaps.  Our  proof here is  similar to the proof of Proposition II.1 in \cite{M}.    

Suppose that  $U_1$ is a gap of $P_1$ and  $U_2$ is a gap of $P'_1$ such $\partial U_1\cap P'_1=\emptyset$ and  $\partial U_2\cap P_1=\emptyset$.
There is two possible case  for a  such pair $(U_1,U_2)$ of linked gaps:
\begin{itemize}
\item   [$a)$]  if  the gap is as the figure below:
 \begin{figure}[ht]
 	\centering % centering figure
 	\scalebox{0.3} % rescale the figure by a factor of 0.8
 	{\includegraphics {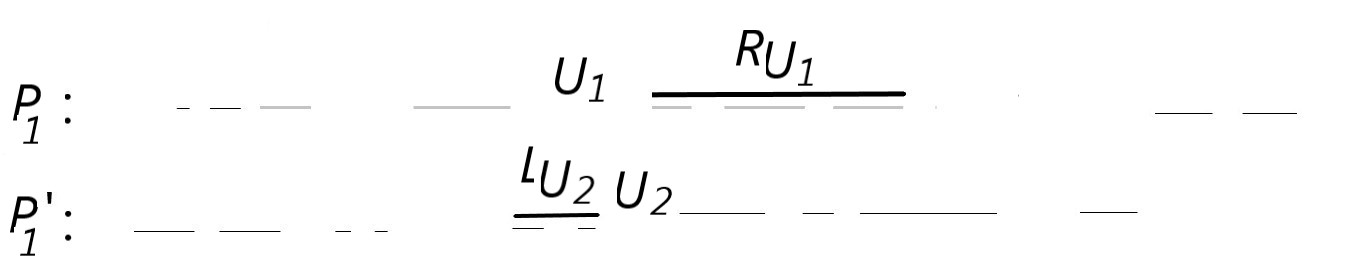}}\caption{\small{For these  gaps, $R(U_1)>|U_2|$ but  $L(U_2)<|U_1|$.}}
 	\label{fig:exm} % labeling to refer it inside the text
 \end{figure} 
 we always  have $R(U_1)>|U_2|$ or  $L(U_2)>|U_1|$, since   $\tau_R(U_1)\tau_L(U_2)>1$ and  this gives a new pair    of linked gaps  with the smaller size, 
\item  [$b)$]   if  the gap is as the figure below:
\begin{figure}[ht]
 	\centering % centering figure
 	\scalebox{0.3} % rescale the figure by a factor of 0.8
 	{\includegraphics {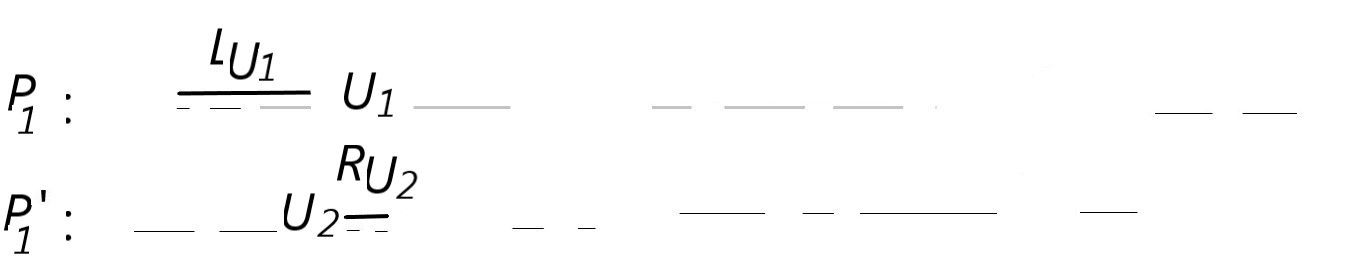}}\caption{\small{For these  gaps, $L(U_1)>|U_2|$ but  $R(U_2)<|U_1|$.}}
 	\label{fig:exm} % labeling to refer it inside the text
 \end{figure} 

 we always  have  $R(U_2)>|U_1|$ or  $L(U_1)>|U_2|$, since   $\tau_L(U_1)\tau_R(U_2)>1$ and  this also  gives a new pair of linked gaps  with the smaller size,
\end{itemize}
In any case, we obtain a smaller pair of gaps, in the sense that the sum of the lengths of its gaps decrease. Using repeatedly the same argument we obtain a  linked gap  $(U_1,U_2)$ such that $\partial U_1\cap P'_1\neq\emptyset$ or $\partial U_2\cap P_1\neq\emptyset$. Thus, $P_1\cap P_1'\neq \emptyset$ which  this says  $P\cap (s P'+t)$ has nonempty interior. This completes the proof of (i).

For $(ii)$,  let  $R_\delta:=\bigcup_{(i,j)\in\Sigma\times\Sigma'}(B^{\delta}_i\cup C^{\delta}_j)$, where $0<\delta$ small and 
\begin{itemize}
\item  $B^{\delta}_i:= \big\{ (s,t)\big|  ~\frac{G(K)}{bp_i}\leq s  \leq \frac{a}{G(K')},~~a_i-bs+\frac{\delta}{p_i} \leq t\leq b_i-\frac{\delta}{p_i}\big\},$
\item $C^{\delta}_j:=\big\{ (s,t)\big|  ~\frac{ G(K)}{b}  \leq  s \leq \frac{aq'_j}{G(K')},~~ -b'_js+\delta\leq t\leq a-a'_js-\delta \big\}$.
\end{itemize}
Obviously, 
 $R$   is a compact set  which  we  will  show that it is a recurrent set, see Figure 3. 
 \begin{figure}[ht]
 	\centering % centering figure
 	\scalebox{0.6} % rescale the figure by a factor of 0.8
 	{\includegraphics {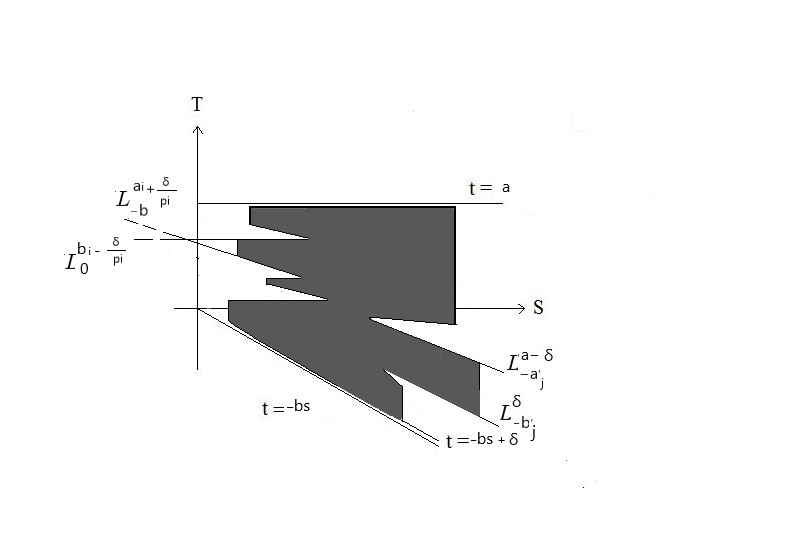}}\caption{\small{Gray region illustrates the recurrent  set $R$.}}
 	\label{fig:exm} % labeling to refer it inside the text
 \end{figure} 
  For given real numbers $c,d\in\mathbb{R},$  let $L_{c}^d:= \{ (s,t)| ~ t=cs+d\}$. Clearly,
\begin{align}
T_i\big(\{L_{c}^d\}\big)=L_{c}^{p_id-p_ia_i},~~~~~~~~~~~~T'_j\big(\{L_{c}^d\}\big)=L_{q_jc+q_ja'j}^{d}.
\end{align} This gives
$$T_i\big(\{L_{-b}^{a_i+\frac{\delta}{p_i} }\}\big)=L_{-b}^{\delta},~~~~~~~~~~~~~T_i\big(\{L_{ 0}^{b_i-\frac{\delta}{p_i}}\}\big)=L^{a-\delta}_{0}$$  and $$T'_j\big(\{L_{ -b'_j}^{\delta}\}\big)=L_{-b}^{\delta},~~~~~~~~~~~~~T'_j\big(\{L_{-a'_j}^{a-\delta}\}\big)=L_{0}^{a-\delta}.$$
Consequently $T_i(B^{\delta}_i)=\{(s,t)| ~s>0,~-bs+\delta\leq t\leq a-\delta\}=T'_j(C^{\delta}_j)$, for all $(i,j)\in\Sigma\times\Sigma'$.

 Since  $(K,K')$ satisfies  the  SGTT and $R$ is compact, one can find  a    small $\delta$ such that for each  $(s,t)\in R$, there is a $(i,j)\in\Sigma\times\Sigma'$
such that $$t\in(a_i+\frac{\delta }{p_i},~b_{i}-\frac{\delta }{p_i})-s (a'_j,~b'_j)=(a_i-b'_js+\frac{\delta }{p_i},~b_i-a'_js-\frac{\delta }{p_i}).$$ This tell that for the such $(s,t)\in R$ and  $(i,j)\in\Sigma\times\Sigma'$
\begin{align}T'_j\circ T_i(s,t)=(\frac{p_is}{q_j},p_it-p_ia_i+p_ia'_js)\in \{\frac{p_is}{q_j}\}\times (-\frac{p_ibs}{q_j}+\delta,~a-\delta).\end{align}

 From   (4.2)  and  the relation $ \frac{q_jG(K)}{bp_i}<\frac{aq_j}{p_iG(K')}$, we obtain that 
 \begin{itemize}
\item  If $\frac{G(K)}{bp_i}\leq s<\frac{aq_j}{p_iG(K')}$, then $T_i(s,t)\in {C^{\delta}_j}^\circ$, since $\pi_1\big(T'_j\circ T_i(s,t)\big)=\frac{p_is}{q_j}< \frac{a}{G(K')}$,
\item If $ \frac{q_jG(K)}{bp_i}< s\leq\frac{aq_j}{G(K')}$,  then $T'_j(s,t){\in B^{\delta}_i}^\circ$, since $\pi_1\big(T_i\circ T'_j(s,t)\big)=\frac{p_is}{q_j}>\frac{G(K)}b{}.$
\end{itemize}
This completes the proof of (ii).

\end{proof}
Thus, if  $(K,K')$ satisfies the   SGTT, then $K$ and  $s K'+t$ have stable intersection,  for all $(s,t)\in R_0^\circ$, where  $R_0=\bigcup_{(i,j)\in\Sigma\times\Sigma'}(B_i\cup C_j)$.  
Theorem \ref{D} introduces  a family of affine Cantor sets  $K$ and $K'$    that may not satisfy GTT,  are linked  together and    do not intersect, while there are translations of  them which  have   stable intersection.
\begin{theorem}\label{D}    Suppose  that  for a   small  $\delta>0$, there are two  open  subsets  $B\subset R_\delta\setminus \cup_{i\in\Sigma}C_i$ and $C\subset R_\delta\setminus\cup_{j\in\Sigma'}B_j$which satisfy 
  \begin{itemize}
\item  [$i)$]  for each $(s,t)\in R_\delta\setminus B \cup C $, the set $P\cap (s P'+t)$ has nonempty interior,
\item  [$ii)$] for each $(s,t)\in \overline B$ and $j\in\Sigma'$, there is $k\in \Sigma$, such that 
$\big(\frac{q_j}{q_k}s,t+(a'_k-a'_j)q_js\big)$ belongs to interior of $R_\delta\setminus B \cup C$, 
\item [$iii)$]  for each $(s,t)\in \overline C$ and $i\in\Sigma$, there is $k\in \Sigma'$, such that  $\big(\frac{p_k}{p_i}s,\frac{p_k}{p_i}t+(a_i-a_k)p_k\big)$  belongs to interior of $ R_\delta\setminus B \cup C$.  
\end{itemize}
 Then  $ R_\delta\setminus B \cup C $ is a recurrent set for  $K$  and $K'$. Thus,   for each $(s,t)\in R_\delta\setminus B \cup C $,  Cantor sets $K$ and $ s K'+t$ have  stable intersection.
\end{theorem} 
\begin{proof} Let   $(s,t)\in R_\delta\setminus B \cup C $, the same as    the proof  stated in  Theorem 3(II), the property $(I)$ gives  a  $(i,j)\in\Sigma\times\Sigma$ such that  
\begin{itemize}
\item  If $\frac{G(K)}{bp_i}\leq s<\frac{aq_j}{p_iG(K')}$, then    $(s_1,t_1)=T_i(s,t)\in {C^{\delta}_j}^\circ$. If $(s_1,t_1)\in \overline C$, then  from $(III)$, $T_k(s,t)\in {\big(R_\delta\setminus B \cup C\big)}^\circ$, since $$T_k\circ {T_i}^{-1}(s_1,t_1)=  \big(\frac{p_k}{p_i}s,\frac{p_k}{p_i}t+(a_i-a_k)p_k\big).$$ 
\item If $ \frac{q_jG(K)}{bp_i}< s\leq\frac{aq_j}{G(K')}$,  then $(s_1,t_1)=T'_j(s,t){\in B^{\delta}_i}^\circ$, If $(s_1,t_1)\in \overline B$, then  from $(II)$, $T'_k(s,t)\in {\big(R_\delta\setminus B \cup C\big)}^\circ$, since
$$T'_k\circ {T'_j}^{-1}(s_1,t_1)=\big(\frac{q_j}{q_k}s,t+(a'_k-a'_j)q_js\big).$$ 
\end{itemize}

These two fact  say that $ R_\delta\setminus B \cup C $ is a recurrent set and completes the proof.
\end{proof}

\begin{proposition} \label{E}  Assume  that  $K$ and $ K'$  are two affine  Cantor sets. Then 
if  $(K,K')$ satisfies the   GTT, then   $K-K'$ is  a   finite union of closed intervals.
\end{proposition} 
\begin{proof} Tending $\delta\to 0$ for the set $R=\bigcup_{(i,j)\in\Sigma\times\Sigma'}(B^{\delta}_i\cup C^{\delta}_j)$   appeared  in $(ii)$,  we obtain  that every  pair $(s,t)$ of set 
  $R_0:=\bigcup_{(i,j)\in\Sigma\times\Sigma'}(B_i\cup C_j)$, where
\begin{itemize}
\item  $B_i:= \big\{ (s,t)\big|  ~\frac{G(K)}{bp_i}\leq s  \leq \frac{a}{G(K')},~~a_i-bs\leq t\leq b_i\big\},$
\item $C_j:=\big\{ (s,t)\big|  ~\frac{ G(K)}{b}  \leq  s \leq \frac{bq'_j}{G(K')},~~ -b'_js\leq t\leq a -a'_js\big\}$.
\end{itemize}
is a difference pair for  $K$ and $K'$ and so $K-sK'=[-bs,~a]$, if and only if $\frac{G(K)}{b}\leq s\leq\frac{a}{G(K')}$.  Also,
  if $s<\frac{G(K)}{b}$   and $(s,t)$ is a difference pair for  $K$ and $K'$ such that $t\in K-sK'$,  then there is  a  finite sequence  $\{a_{i}\}^{m}_{i=1}$  of the set $\Sigma$ and a  ${(i,j)\in\Sigma\times\Sigma'}$ such that 
$$T_{a_{m}}\circ T _{a_{m-1}} \circ\ldots\circ
T_{a_{1}}(s,t)\in B_i\cup C_j.$$
Let $p=\max_{i\in \Sigma} \{p_i\}$ and  $q=\max_{j\in \Sigma'} \{q_j\}$, 

$$P=\max\{\Pi_{i_j\in\Sigma} p_{i_1}p_{i_2}\dots p_{i_k}|~p<p_{i_1}p_{i_2}\dots p_{i_k}s<q\}$$

Thus, the connected component of the set $K-sK'$ containing  $t$ have at least the size $\frac{1}{P}\big(\frac{a}{p}+G(K)\big)$ (see Definition of $B_i$) and so   $K-sK'$ is  a   finite union of closed intervals.  A similar proof can be written for the case $s> \frac{1}{G(K')}$  and this completes the proof of theorem.
\end{proof}

\end{document}